\documentclass[12pt,a4paper]{article}

\usepackage{graphicx}
\usepackage{amssymb}
\usepackage{amsmath}
\usepackage{epstopdf}

\usepackage{float}

\usepackage{caption}
\usepackage{subcaption}

\usepackage[margin=0.8in]{geometry} 

\usepackage{subfloat}

\title{The Unified Transform Method: beyond circular or convex domains}

\author{Jesse J. Hulse$^{1}$, Loredana Lanzani$^{2}$, Stefan G. Llewellyn Smith$^{3,4}$ \& Elena Luca$^{5}$}

\date{}

\graphicspath{{Fig/}}

\numberwithin{equation}{section}

\usepackage{pst-plot,booktabs,mathtools}
\usepackage{amsfonts, amsthm, xcolor}
\usepackage{graphicx}
\newcommand{\A}{\mathcal{A}}
\newcommand{\B}{\mathcal{B}}

\newcommand{\dd}{\mathrm{d}}
\newcommand{\ee}{\mathrm{e}}
\newcommand{\ci}{\mathrm{i}}
\DeclareMathOperator*{\Arg}{arg}

\DeclareMathOperator*{\Imag}{Im}
\DeclareMathOperator*{\Real}{Re}

\newtheorem{lemma}{Lemma}

\newtheorem{theorem}{Theorem}

\usepackage{color}
\newcommand{\Blue}[1]{{\color{blue}{#1}}}
\newcommand{\Red}[1]{{\color{red}{#1}}}

\usepackage{subcaption}

\begin{document}

\maketitle


\begin{center}
$^{1}$Department of Mathematics \\
Syracuse University \\
Syracuse, NY 13244-1150, USA
\end{center}

\begin{center}
$^{2}$Department of Mathematics \\
University of Bologna \\
Bologna, Italy
\end{center}

\begin{center}
$^{3}$Department of Mechanical and Aerospace Engineering \\
Jacobs School of Engineering, UCSD \\
La Jolla, CA 92093-0411, USA.
\end{center}

\begin{center}
$^{4}$Scripps Institution of Oceanography, UCSD \\
La Jolla, CA 92039-0213, USA.
\end{center}

\begin{center}
$^{5}$Climate and Atmosphere Research Center \\
The Cyprus Institute \\
Nicosia, 2121, Cyprus \\
\vskip 0.05truein 
Corresponding author: {\tt e.louca@cyi.ac.cy}
\end{center}

\vskip 0.5truein
\begin{center}
{\bf Subject Areas}
\end{center}
\begin{center}
  applied mathematics, transform
 methods, analysis
\end{center}

\begin{center}
{\bf Keywords}
\end{center}
\begin{center}
 complex analysis, unified transform
 method, Cauchy kernel, Szeg\H o kernel
\end{center}

\begin{center}
{\bf Abstract}
\end{center}
A new transform-based approach is presented that can be used to solve mixed boundary value problems for Laplace's equation in non-convex and other planar domains, specifically the so-called Lipschitz domains. This work complements Crowdy (2015, CMFT, 15, 655--687), where new transform-based techniques were developed for boundary value problems for Laplace's equation in circular domains. The key ingredient of the present method is the exploitation of the properties of the Szeg\H o kernel and its connection with the Cauchy kernel to obtain transform pairs for analytic functions in such domains. Several examples are solved in detail and are numerically implemented to illustrate the application of the new transform pairs.
\noindent

\vspace{1cm}
\noindent

\vfill\eject

\section{Introduction}
The Unified Transform Method (UTM), introduced by A.~S. Fokas in the late 1990s \cite{Fokas1997}, is a technique for analyzing boundary value problems for linear and integrable nonlinear PDEs. Since its inception, the UTM has garnered significant interest within the applied mathematics community. Over time, numerous adaptations of the original method have been developed to address specific classes of equations. For the Laplace, biharmonic, Helmholtz, and modified Helmholtz equations in convex polygonal domains, the UTM yields integral representations of solutions in the complex Fourier plane \cite{FokasKapaev2003, CrowdyFokas2004,  CrowdyLuca:2014, DimakosFokas2015, SpenceFokas2010, DavisFornberg2014}.

In particular, for Laplace's equation, Fokas \& Kapaev \cite{FokasKapaev2003} formulated a transform method to solve boundary value problems in simply-connected convex polygons. Their approach initially utilized various techniques, including spectral analysis of parameter-dependent ODEs and Riemann--Hilbert methods. Later, Crowdy \cite{Crowdy:2015, Crowdy:2015b} demonstrated that this method could be reformulated using a complex function-theoretic framework, leading to the development of a new transform pair tailored to circular domains (domains bounded by circular arcs, with line segments and hence convex polygons as a special case). The authors of the present paper have recently obtained a new construction of so-called \emph{quasi-pairs} that extends the original approach of Fokas \& Kapaev \cite{FokasKapaev2003} from convex {\it polygons} to {\it any} convex domains \cite{HLLL:2024}, e.g., to ellipses.

The focus of the current study is the development of a transform-based method for Laplace's equation in planar domains that may fail to be convex or circular, specifically the class of so-called Lipschitz domains \cite{Pommerenke:1992}. Hence this new method provides means for solving mixed boundary value problems in non-convex domains, such as in the interior of a concave quadrilateral. To achieve this, we build upon Crowdy's \cite{Crowdy:2015} construction for circular domains and develop a new transform method for solving mixed boundary value problems in Lipschitz domains. We obtain tailor-made transform pairs to represent the value of a given function at any point in the domain of interest by exploiting the properties of the Szeg\H o kernel for the domain in question, along with an appropriate conformal mapping.

The Szeg\H o and Cauchy kernels are important in complex analysis and operator theory. Among their many features, they can be used to produce and reproduce holomorphic i.e., analytic functions in a given domain $D$ from the functions' known boundary values (the classical Cauchy integral formula is one instance). Both kernels are associated to projections from the Hilbert space of square integrable functions $L^2(\partial D)$ to the Hardy space $H^2(D)$ \cite{HLLL:2024}. (Here $\partial D$ denotes the boundary of $D$.) However, while the Cauchy kernel is completely explicit and canonical (it does not change with $D$), the Szeg\H o kernel is domain-specific and in general it is not known explicitly. On the other hand the Cauchy kernel does not have a good transformation law under conformal maps, $\Phi : D \rightarrow \mathbb{D}$ (where $\mathbb{D}$ is the unit disc), but the Szeg\H o kernel does. It turns out that the Szeg\H o and Cauchy kernels associated to $D$ are identical precisely when $D$ is the unit disc, i.e.~$D=\mathbb{D}$ \cite{KerzmanStein:1978}. This means that whenever  one has enough knowledge of the conformal map for the given domain $D$, one can use it to obtain an accurate numerical approximation of the Szeg\H o kernel for $D$. This idea makes it possible to extend the scope of the UTM devised by Crowdy \cite{Crowdy:2015} for circular domains, to any domain $D$ for which one has a satisfactory understanding of its conformal map to the unit disc. In this paper, we develop this idea for a general simply-connected domain $D$ whose boundary satisfies the condition known as ``Lipschitz regularity'' \cite{Pommerenke:1992}, and then we carry it out explicitly for specific choices of $D$. These include examples where $D$ is an ellipse, for which a different approach was developed in our previous paper \cite{HLLL:2024}, as well as the cases where $D$ is a concave polygon or a punctured disc.

In \S\,\ref{background}, we present the theoretical background needed to formulate the new transform pairs for analytic functions in Lipschitz planar domains. The transform pairs are then constructed in \S\,\ref{sectiontransform} for simply-connected, Lipschitz domains and in \S\,\ref{sectiontransformpunctured} for punctured domains. The next step involves implementing the method for a variety of mixed boundary value problems in \S\,\ref{sectionapplication1}--\ref{sectionapplication2}.  Finally, we conclude and discuss further applications in \S\,\ref{sectiondiscussion}. (Some technical details are relegated to Appendices \ref{Appendix A} and \ref{Appendix B}.)

\section{Background \label{background}}\label{sec2}
Let $D \subset \mathbb{C}$ denote a bounded domain in $\mathbb C$ whose boundary is denoted by $\partial D$ and closure by $\overline{D}$, that is $\overline D = D\cup \partial D$. In what follows, interior points are denoted by $z \in D$ and boundary points by $\zeta \in \partial D$. To begin with, we assume that $D$ is a sufficiently smooth domain, say $D$ is of class $C^1$, so that the unit normal and tangent vectors to $\partial D$ are continuously differentiable functions of $\zeta \in \partial D$. The unit tangent vector at any point $\zeta \in \partial D$ is denoted by $T_D(\zeta)$ and is given by
\begin{equation}
T_D(\zeta)=\frac{\zeta'(t)}{|\zeta'(t)|},
\label{2.1}
\end{equation}
where $\zeta=\zeta(t)$ denotes any parametrization for $\partial D$. Finally, the arc-length measure for $\partial D$ is denoted by $\sigma$, and hence
\begin{equation}
\dd \sigma(\zeta)=|\zeta'(t)|\dd t.
\label{2.2}
\end{equation}
It is important to point out that the above and what we say below can also be made meaningful for less smooth domains, such as the Lipschitz domains \cite{Pommerenke:1992} in which case, e.g., \eqref{2.1} and \eqref{2.2} will be valid for $\sigma$-almost every $\zeta \in \partial D$. See also the note on the regularity of the conformal mapping in \S\,\ref{Szego as alternative}\ref{advantages of Szego} below.

Before we proceed further, a couple of definitions are introduced. (See \cite{Duren:1970}, Chapter 10, for a more detailed discussion.) Let $D$ be a simply-connected domain with at least two boundary points and let $0<p\leq \infty$. An analytic function $f$ on $D$ belongs to $H^p(D)$ if the subharmonic function $|f|^p$ has a harmonic majorant in $D$. The function $f$ belongs to $E^p(D)$ if there exists a sequence of rectifiable Jordan Curves $C_1,C_2,\dots$ in $D$ converging to the boundary in the sense that $C_n$ eventually contains each compact subdomain of $D$, such that
\begin{equation}
\sup_{n} \int_{C_n}|f(z)|^p \, |\dd z|\leq M<\infty.
\label{Hp condition}
\end{equation}
It is a well-known fact that $E^p(D)$ and $H^p(D)$ coincide if $D$ is the unit disk, and the same holds true for smooth domains \cite{Duren:1970}. On the other hand, one has that $E^{p}(D) \neq H^{p}(D)$ for Lipschitz $D$ (see \cite{Duren:1970}, page 183) and we henceforth work with $E^{p}(D)$. If $q \geq p$ then $E^q(D) \subseteq E^p(D)$. In particular we have that $E^p(D) \subseteq E^1(D)$ for any $p \geq 1$; moreover, analytic functions that are continuous on $\overline{D}$ belong to $E^p(D)$ for any $p>0$.

\subsection{The Cauchy integral formula}

\subsubsection{Definition}

For a function $f \in E^1(D)$, Cauchy's integral formula states that, for $z \in D$,
\begin{equation}\label{E:CF}
f(z)=\frac{1}{2\pi\ci} \int_{\partial D} {\frac{f(\zeta)}{\zeta-z} \,\dd \zeta}.
\end{equation}

\subsubsection{The Cauchy integral formula as an inner product}

In the sequel it will be convenient to express the Cauchy integral formula \eqref{E:CF} in terms of the inner product in  $L^2 (\partial D, \dd \sigma)$. To this end, recall that the quantity $\dd \zeta$ in \eqref{E:CF}  is the complex differential $\dd \zeta = \zeta'(t) \, \dd t$, hence
\begin{equation}\label{E:dz-ds}
\overline{T_D(\zeta)} \,\dd \zeta  = \dd \sigma (\zeta), \quad \text{or, equivalently,} \quad \dd \zeta =  T_D(\zeta) \, \dd \sigma (\zeta),
\end{equation}
where $\overline{a}$ denotes the complex conjugate of the quantity $a \in \mathbb{C}$. 
Using \eqref{E:dz-ds}, Cauchy's integral formula \eqref{E:CF} can be rewritten as
\begin{equation}
f(z)=\int_{\partial D} {f(\zeta) \frac{T_D(\zeta)}{2\pi\ci(\zeta-z)}\,\dd \sigma(\zeta)}=\int_{\partial D} {f(\zeta) \overline{ \left[ \frac{\overline{T_D(\zeta)}}{-2\pi\ci(\overline{\zeta}-\overline{z})} \right] } \,\dd \sigma(\zeta)}.
\end{equation}
Thus Cauchy's integral formula \eqref{E:CF} may be reinterpreted as follows:

\begin{theorem}{(Cauchy formula)}
For any $f \in E^1(D)$, we have
\begin{equation}
f(z)=\int_{\partial D} f(\zeta) \overline{C_{D}(\zeta, z)} \,\dd \sigma(\zeta) =: \langle f, C_{D}(\cdot, z) \rangle_{L^2(\partial D, \dd \sigma)}, \quad z \in D,
\label{Cauchyformula1}
\end{equation}
where $\langle \cdot, \cdot \rangle_{L^2(\partial D, \dd \sigma)}$ denotes the inner product in $L^2(\partial D, \dd \sigma)$, and $C_{D}(\zeta, z)$ is the {\bf Cauchy kernel} for $D$ \cite{Bell:1992} and is given by
\begin{equation}
C_{D}(\zeta, z) \coloneqq -\frac{\overline{T_{D}(\zeta)}}{2\pi\ci(\overline{\zeta}-\overline{z})}, \quad \zeta \in \partial D, \quad z \in D.
\label{Cauchykernel1}
\end{equation}
\end{theorem}
The reason for this inner product-based interpretation of the Cauchy formula and kernel \eqref{Cauchyformula1}--\eqref{Cauchykernel1}, is because the Szeg\H o formula and kernel in the next subsection can only be defined in terms of inner product in $L^2(\partial D, \dd \sigma)$. Thus \eqref{Cauchykernel1} will be essential in order to obtain an explicit formula for the Szeg\H o kernel. See \eqref{E:Szego-disc} and \eqref{Szego2.17} below.

\subsection{The Szeg\H o kernel as an alternative for the Cauchy kernel \label{Szego as alternative}}

Let $E^2(D)$ be the Hardy space previously defined.

\subsubsection{Definition and examples}

\begin{theorem}{(Szeg\H o formula)}
For any $f$ in the Hardy space $E^2(D)$, we have
\begin{equation}\label{E:SF}
f(z)=\int_{\partial D} {f(\zeta) \overline{S_{D}(\zeta, z)} \,\dd \sigma(\zeta)} = \langle f, S_{D}(\cdot, z) \rangle_{L^2(\partial D, \dd \sigma)}, \quad z \in D,
\end{equation}
where $S_{D}(\zeta, z)$ is the {\bf Szeg\H o kernel for $D$} \cite{Bell:1992}.
\end{theorem}

\noindent
We point out that, while the definition of the Cauchy kernel is completely explicit for any domain $D$ (see \eqref{Cauchykernel1}), there is no explicit formula that directly provides the Szeg\H o kernel for an arbitrary domain $D$. More precisely, the {\it existence} of the Szeg\H o kernel is established by the Riesz representation theorem \cite{Conway:1990} which is only an existence theorem and does not provide a formula-based definition for $S_D(\zeta,z)$; see \cite{Bell:1992} for a complete treatment of the Szeg\H o kernel.

There are a few specific domains whose Szeg\H o kernel can be computed explicitly, a number of examples are listed below:

\begin{itemize}
\item[(a)] It was proved by Kerzman \& Stein \cite{KerzmanStein:1978} that the only domain for which the Szeg\H o kernel agrees with the Cauchy kernel is the unit disc. To be precise, if $D =\mathbb D \equiv \{z \in \mathbb{C} | |z|<1\}$, then
\begin{equation}\label{E:Szego-disc}
S_{\mathbb{D}}(\zeta, z)=C_{\mathbb{D}}(\zeta, z)=-\frac{\overline{T_\mathbb{D}(\zeta)}}{2\pi\ci(\overline{\zeta}-\overline{z})} = \frac{1}{2\pi}\frac{1}{(1-\zeta \overline{z})},\quad |\zeta|=1, \quad |z|<1,
\end{equation}
where we have used that $T_\mathbb{D}(\zeta)=\ci \zeta$. Hence, if $D=\mathbb{D}$, then the Szeg\H o kernel for $\mathbb D$ equals the Cauchy kernel.
\end{itemize} 

\noindent
It turns out that the Hardy space $E^p(D)$ and the Szeg\H o kernel may be defined also for so-called punctured domains \cite{GallagherGuptaLanzaniVivas:2021}. We recall known examples in (b) and (c) below.

\begin{itemize}
\item[(b)] If $D= \mathbb{D}^{*} \equiv \{z \in \mathbb{C} | 0<|z|<1\}$ (the punctured disc) and $f$ has a pole of order $n\in\mathbb N$ at $0$ and is otherwise analytic on $\mathbb D^*$ and  continuous on $\overline{\mathbb{D}^{*}}$, then \eqref{E:SF} is valid when interpreted as an integral on $\partial \mathbb{\mathbb{D}}$ with
\begin{equation}\label{E:Szego-punctured-disc}
S_{\mathbb{D^*}}(\zeta, z) \equiv S^{(n)}_{\mathbb{D^*}}(\zeta, z) \coloneqq 
\frac{1}{2\pi}\, \frac{1}{(\zeta \overline{z})^n(1-\zeta\overline{z})}, \quad |\zeta|=1, \quad 0<|z|<1.
\end{equation}
That is, for any $f$ which has a pole of order $n$ at $0$, we have
\begin{equation}
f(z)=\int_{\partial \mathbb{D}} {f(\zeta) \overline{S^{(n)}_{\mathbb{D}^{*}}(\zeta, z)} \,\dd \sigma(\zeta)} =
\langle f, S^{(n)}_{\mathbb{D}^*}(\cdot, z) \rangle_{L^2(\partial \mathbb{D}, \dd \sigma)}, \quad z \in \mathbb{D}^{*}.
\end{equation}

\item[(c)] More generally, if $D$ is a simply-connected domain and $D^*=D\setminus\{z_0\}$, for $z_0\in D$, is the punctured domain, and $f \in E^1(D^{*})$ has a pole of order $n$ at $z_0$, then
\begin{equation}
f(z)=\int_{\partial D} {f(\zeta) \overline{S^{(n)}_{D^{*}}(\zeta, z)} \,\dd \sigma(\zeta)} =
\langle f, S^{(n)}_{D^*}(\cdot, z) \rangle_{L^2(\partial D, \dd \sigma)}, \quad z \in D^{*},
\end{equation}
with
\begin{equation}\label{E:Szego-punctured}
S^{(n)}_{{D^*}}(\zeta, z) = \left(M_{z_0}\circ\Phi\right)^{-n} (z)\ S_D(\zeta, z)\, \overline{\left(M_{z_0}\circ\Phi\right)^{-n} (\zeta)}, \quad \zeta \in \partial D, \quad z \in D^{*},
\end{equation}
where $S_D(\zeta, z)$ is the Szeg\H o kernel for $D$, $\Phi: D\to\mathbb D$ is a conformal mapping, and
\begin{equation}
M_{z_0}(z)\coloneqq \frac{z-z_0}{1-\overline{z_0}z}.
\end{equation}
The expression \eqref{E:Szego-punctured-disc} is a special case of \eqref{E:Szego-punctured}, with $D=\mathbb D$, $z_0=0$ and $\Phi(z)=z=M_{z_0}(z)$. Also, note that an analogous formula holds for $D^*=D\setminus\{z_1,\ldots, z_m\}$ (a simply-connected domain with $m$ punctures) \cite{GallagherGuptaLanzaniVivas:2021}.

\end{itemize}

\subsubsection{Advantages of the Szeg\H o kernel}
The advantage of  the Szeg\H o kernel is that it enjoys a good transformation law under conformal maps between domains in the complex plane, whereas the Cauchy kernel does not. Namely, if $\Phi: D_1\mapsto D_2$ is conformal, then
 \begin{equation}\label{E:Szegokernel}
S_{D_1}(\zeta,z)\coloneqq \overline{\sqrt{\Phi'(z)}} S_{D_2}(\Phi(\zeta),\Phi(z)) \sqrt{\Phi'(\zeta)}, \quad \zeta \in \partial D_1, \quad z \in D_1,
\end{equation} 
which is valid, for example,
\begin{itemize}
\item[(a)] for any $z \in D_1$ and any $\zeta \in \partial D_1$, if $D_1$ is simply-connected and of class $C^{1,1}$, and $D_2=\mathbb{D}$ (the classical setting);

\item[(b)] if $D_1$ is simply-connected and Lipschitz, and $D_2=\mathbb{D}$, then equation \eqref{E:Szegokernel} is valid in the sense that, for any $z\in D_1$,
\begin{equation}
S_{D_1}(\cdot, z)\coloneqq \overline{\sqrt{\Phi'(z)}} S_{\mathbb{D}}(\Phi(\cdot),\Phi(z)) \sqrt{\Phi'(\cdot)} \quad \text{as functions in }
L^2(\partial D_1, \dd \sigma).
\label{Szego2.17}
\end{equation}
It follows that, for any $z \in D_1$, we have
\begin{equation}\label{E:Szegokernel-Lipschitz}
S_{D_1}(\zeta, z) \equiv \overline{\sqrt{\Phi'(z)}} S_{\mathbb{D}}(\Phi(\zeta),\Phi(z)) \sqrt{\Phi'(\zeta)}, \quad \text{for }\ \sigma \text{-almost every } \zeta \in \partial D_1,
\end{equation}
and this formula for $S_{D_1}(\zeta,z)$ can be made even more explicit by employing the expression \eqref{E:Szego-disc} for $S_{\mathbb{D}}(\Phi(\zeta),\Phi(z))$.

\item[(c)] if $D_1 = \mathbb A_{1,\infty}(0) \equiv \{z \in \mathbb{C} | 1<|z|<\infty\}$ and $D_2=\mathbb{D}^*\coloneqq \mathbb{D} \setminus \{0\}$, then \eqref{E:Szegokernel}  is valid for any $\zeta \in \partial \mathbb A_{1,\infty}(0)$ and any $z \in \mathbb A_{1,\infty}(0)$ with $\Phi(z) \equiv 1/z$; see \cite{TegtmeyerThomas:1999}.

\end{itemize}

\noindent
Returning to \eqref{E:Szegokernel}, using \eqref{E:dz-ds}, we can write
\begin{equation}\label{E:Szego-cont}
\overline{S_{D_1}(\zeta,z)}\,\dd \sigma(\zeta) =\sqrt{\Phi'(z)} ~\overline{S_{D_2}(\Phi(\zeta),\Phi(z))}~\overline{{\sqrt{\Phi'(\zeta)}}}~\overline{T_{D_1}(\zeta)}\, \dd \zeta.
\end{equation}
Furthermore it is known (\cite{Bell:1992}, p. ~53) that
\begin{equation}\label{E:tg-vecs}
T_{ D_2}(\Phi (\zeta))\overline{\sqrt{\Phi'(\zeta)}}  =
 \sqrt{\Phi'(\zeta)}T_{D_1}(\zeta),\quad \zeta \in \partial D.
\end{equation}
Substitution of \eqref{E:tg-vecs} into \eqref{E:Szego-cont} gives
\begin{equation}\label{E:Szego-final}
\overline{S_{D_1}(\zeta,z)}\,\dd \sigma(\zeta) =\sqrt{\Phi'(z)}~\overline{S_{D_2}(\Phi(\zeta),\Phi(z))}~{\sqrt{\Phi'(\zeta)}}~\overline{T_{D_2}(\Phi(\zeta))} \dd \zeta.
\end{equation}

\hspace{0.1cm}
\subsubsection{Regularity of conformal mapping \label{advantages of Szego}}
The results stated in the following two bullets will be needed in order to justify the validity of e.g., \eqref{3.11}--\eqref{statement} in the proof of Lemma \ref{lemma 3.1} below.

\begin{itemize}
\item Pommerenke \cite{Pommerenke:1992} (Theorems 3.5 and 3.6): Let $D\subset\mathbb C$ be a proper simply-connected domain of class $C^{1,1}$, and let $\Phi: D\to\mathbb D$ be a conformal map. Then $\Phi$ extends to a map: $\overline{D}\to \overline{\mathbb D}$ with $\Phi: \partial D\to \partial \mathbb D$. Furthermore, $\Phi \in C^1(\overline D)$ and $\Phi'(\zeta)\neq 0$ for any $\zeta\in \overline D$.
\item Duren \cite{Duren:1970} (Theorem 3.12): Let $D\subset\mathbb C$
be a proper simply-connected domain, let $\Phi: D\to\mathbb D$ be a conformal map, and let $\Psi: \mathbb D\to D$ be the inverse map. Then $\partial D$ is rectifiable if and only if $\Psi'\in H^1(\partial \mathbb D)$. Cauchy's theorem holds on domains that have rectifiable boundary (Duren \cite{Duren:1970}; Theorem 10.4). Lipschitz domains have rectifiable boundary, thus Cauchy's theorem may be applied in Lemma \ref{lemma 3.1} below. 
\end{itemize}

\subsection{Transform pair for circular domains via the Cauchy kernel function}

\subsubsection{The transform pair}

Returning to the Cauchy integral formula \eqref{E:CF}, for the case when $D=\mathbb{D}$ is the unit disc, it was proved by Crowdy \cite{Crowdy:2015} that the expression $1/(\zeta-z)$ (for $\zeta \in \partial \mathbb{D}$ and $z \in \mathbb{D}$; see Fig.~\ref{geometric construction}) admits the spectral representation
\begin{equation}\label{E:spectral-disc}
\frac{1}{\zeta-z}=\int_{L_{1}} {\frac{1}{1-\ee^{2\pi \ci k}}\, \frac{z^{k}}{\zeta^{k+1}} \,\dd k}+\int_{L_{2}} {\frac{z^{k}}{\zeta^{k+1}} \,\dd k}+\int_{L_{3}} {\frac{\ee^{2\pi \ci k}}{1-\ee^{2\pi \ci k}}\, \frac{z^{k}}{\zeta^{k+1}} \,\dd k},
\end{equation}
where $L_j$, $j=1$, $2$, $3$ are contours in the spectral $k$-plane, as shown in Fig.~\ref{fundamental contour}. That is, the contour $L_1$ is the union of the negative imaginary axis $(-\ci \infty, \ci r]$ and the arc of the quarter circle $|k|=r$, $0<r<1$, in the third quadrant traversed in a clockwise sense; the contour $L_2$ is the real interval $[-r,\infty)$; the contour $L_3$ is the arc of the quarter circle $|k|=r$ in the second quadrant traversed in a clockwise sense together with the portion of the
positive imaginary axis $[\ci r, \ci \infty)$ \cite{Crowdy:2015}; see Fig.~ \ref{fundamental contour}. Replacing the function $1/(\zeta-z)$ in \eqref{E:CF} with its spectral representation \eqref{E:spectral-disc}, we find the following {\bf transform pair for the interior of the unit disc} \cite{Crowdy:2015}:
\begin{equation}
\begin{cases}
f(z)=\displaystyle \frac{1}{2\pi \ci} \left[ \int_{L_{1}}{\frac{\rho(k)}{1-\ee^{2\pi \ci k}}z^{k} \dd k}+\int_{L_{2}} {\rho(k) z^{k} \,\dd k} +\int_{L_{3}}{\frac{\rho(k) \ee^{2\pi \ci k}}{1-\ee^{2\pi \ci k}}z^{k} \,\dd k}\right], \quad z \in \mathbb{D}, \\ \\
\rho(k)\coloneqq \displaystyle \oint_{\partial \mathbb{D}}{\frac{f(\zeta)}{\zeta^{k+1}}\,\dd \zeta}, \quad k \in \mathbb{C}.
\end{cases}\label{Crowdy's original transform pair}
\end{equation}
The {\bf global relation} states that
\begin{equation}
\rho(k)=0, \quad k \in \mathbb{Z}^{-} \equiv \{\cdots,-3,-2,-1\};
\end{equation}
see \cite{Crowdy:2015}.

\begin{figure*}[t!]
    \centering
    \begin{subfigure}[t]{0.5\textwidth}
        \centering
        \includegraphics[width=0.6\textwidth]{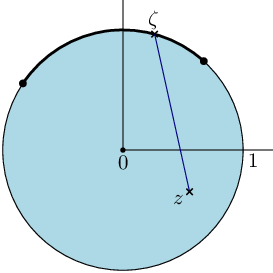}
        \caption{}
        \label{geometric construction}
    \end{subfigure}%
    ~
    \begin{subfigure}[t]{0.5\textwidth}
        \centering
        \includegraphics[width=1\textwidth]{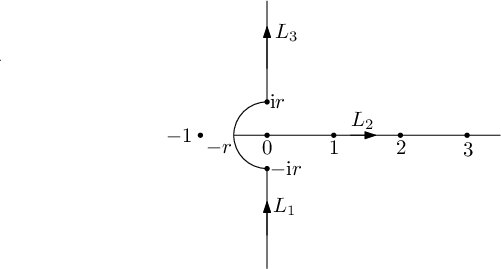}
        \caption{}
        \label{fundamental contour}
    \end{subfigure}
    \caption{(a) Point $\zeta$ on the boundary of $\mathbb{D}$ and point $z$ in its interior, (b) The fundamental contour for circular arc edges with $0<r<1$ \cite{Crowdy:2015}.}
\end{figure*}

\begin{figure*}[t!]
    \centering
    \begin{subfigure}[t]{0.5\textwidth}
        \centering
        \includegraphics[width=0.6\textwidth]{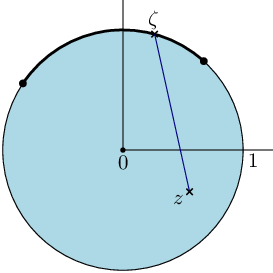}
        \caption{}
        \label{geometric construction}
    \end{subfigure}%
    ~
    \begin{subfigure}[t]{0.5\textwidth}
        \centering
        \includegraphics[width=1\textwidth]{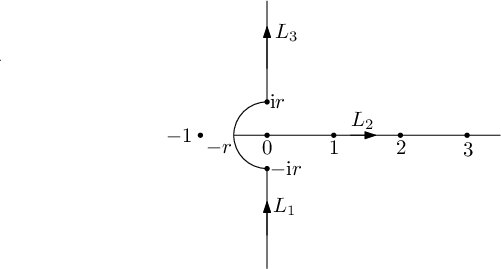}
        \caption{}
        \label{fundamental contour}
    \end{subfigure}
    \caption{(a) Point $\zeta$ on the boundary of $\mathbb{D}$ and point $z$ in its interior, (b) The fundamental contour for circular arc edges with $0<r<1$ \cite{Crowdy:2015}.}
\end{figure*}

 \section{A new transform pair and global relation \label{sectiontransform}}
 
In this section, we show how to obtain transform pairs for simply-connected Lipschitz domains by exploiting the properties of the Szeg\H o kernel.

\subsection{A transform pair for Lipschitz domains}
\begin{theorem}
    Suppose that $D$ is a bounded simply-connected domain with Lipschitz boundary, let $\Phi: D \mapsto \mathbb{D}$ be a conformal map, and let $f: D \mapsto \mathbb C$ be in the Hardy space $E^2(D)$. Then the following holds:
\begin{equation}\label{Transform pair domain1}
f(z)=\displaystyle \frac{1}{2\pi {\rm i}} \sqrt{\Phi'(z)} \left[ \int_{L_{1}}{\frac{\rho(k)}{1-\ee^{2\pi \ci k}}\Phi(z)^{k} \,\dd k}+\int_{L_{2}} {\rho(k) \Phi(z)^{k} \,\dd k} +\int_{L_{3}}{\frac{\rho(k) \ee^{2\pi \ci k}}{1-\ee^{2\pi \ci k}}\Phi(z)^{k} \,\dd k}\right], \quad z \in D,
\end{equation}
where
\begin{equation}\label{Transform pair domain2}
\rho(k)\coloneqq \displaystyle \int_{\partial D} {\frac{\sqrt{\Phi'(\zeta)}}{\Phi(\zeta)^{k+1}}f(\zeta) \,\dd\zeta}, \quad k \in \mathbb{C}.
\end{equation}
We refer to the pair $(f,\rho)$ as the transform pair for $D$.

\end{theorem}

\begin{proof}
Combining the formula \eqref{E:Szego-final} for $D_1=D$ and $D_2=\mathbb D$ with the formula \eqref{E:Szego-disc} for $S_{\mathbb D}(\zeta, z)$, we find
\begin{equation}
\begin{split}
\overline{S_D(\zeta, z)} \dd \sigma(\zeta) &= \sqrt{\Phi'(z)}\frac{1}{2\pi \ci} \frac{T_{\mathbb D}(\Phi(\zeta))}{(\Phi(\zeta)-\Phi(z))}\sqrt{\Phi'(\zeta)}\overline{T_{\mathbb D}(\Phi(\zeta))} \dd \zeta \\
 &= \frac{1}{2\pi \ci} \sqrt{\Phi'(z)} \frac{\sqrt{\Phi'(\zeta)}}{(\Phi(\zeta)-\Phi(z))} \dd \zeta.
\end{split}
\end{equation}
Substitution of this expression in the Szeg\H o formula \eqref{E:SF} gives
\begin{equation}\label{E:Szego-conf}
\begin{split}
f(z)&= \int_{\partial D} {f(\zeta) \overline{S_{D}(\zeta,z)} \,\dd \sigma(\zeta)} \\
&=\frac{1}{2\pi \ci} \sqrt{\Phi'(z)} \int_{\partial D} {f(\zeta) \left[ \frac{\sqrt{\Phi'(\zeta)}}{\Phi(\zeta)-\Phi(z)} \right] \,\dd \zeta}, \qquad z \in D.
\end{split}
\end{equation}

Next we point out that $\Phi(\zeta)\in \partial \mathbb D$ and $\Phi (z) \in \mathbb D$ by the mapping properties of $\Phi$, so the spectral decomposition \eqref{E:spectral-disc} gives
\begin{equation}\label{E:spectral-general}
\frac{1}{\Phi(\zeta)-\Phi(z)}=\int_{L_{1}} {\frac{1}{1-\ee^{2\pi{\rm i}k}}\, \frac{\Phi(z)^{k}}{\Phi(\zeta)^{k+1}} \,\dd k}+\int_{L_{2}} {\frac{\Phi(z)^{k}}{\Phi (\zeta)^{k+1}} \dd k}+\int_{L_{3}} {\frac{\ee^{2\pi{\rm i}k}}{1-\ee^{2\pi\ci k}}\, \frac{\Phi (z)^{k}}{\Phi (\zeta)^{k+1}} \,\dd k},
\end{equation}
with $\zeta \in \partial D$, $z \in D$ and where the contours $L_j$, $j=1,\ldots, 3$ are shown in Fig.~\ref{fundamental contour}. Combining \eqref{E:spectral-general} with \eqref{E:Szego-conf}, we obtain \eqref{Transform pair domain1} with $\rho(k)$ as in \eqref{Transform pair domain2}.
\end{proof}

\subsection{A new global relation}
\begin{lemma}
Let $D$ be a bounded Lipschitz domain, $\Phi: D \mapsto \mathbb{D}$ be a conformal map and $f \in E^{2}(D)$. Then,
\begin{equation}
\frac{\sqrt{\Phi'(z)}}{\Phi^k(z)} f(z) \in H^2(D), \quad z \in D, \quad \forall k \in \mathbb{Z}^{-}.
\end{equation}
Moreover, for $\rho(k)$ as in \eqref{Transform pair domain2}, we have
\begin{equation}\label{GR}
\rho(k)=0, \quad k \in \mathbb{Z}^{-}.
\end{equation}
We refer to \eqref{GR} as the global relation for $D$.
\label{lemma 3.1}
\end{lemma}

\begin{proof}
Let $k \in \mathbb{Z}^{-}$ and write $n \coloneqq -k-1 \in \mathbb{Z}^{+}_{0}\equiv \{0,1,2,3,\cdots\}$. We can then write
\begin{equation}
\tau(n)\coloneqq \rho(-n-1) = \int_{\partial D} {\sqrt{\Phi'(\zeta)} f(\zeta) \Phi(\zeta)^n \,\dd \zeta}.
\end{equation}
In order to prove that $\rho(k)$ defined as in \eqref{Transform pair domain2} satisfies the global relation \eqref{GR}, it is enough to show that
\begin{equation}
\tau(n)=0, \qquad \text{for } n \in \mathbb{Z}^{+}_0.
\label{GR tau(n)}
\end{equation}
We denote the inverse of the conformal map $\Phi$ by $\Psi:=\Phi^{-1}:\mathbb{D} \rightarrow D$. By Caratheodory's theorem \cite{Duren:1970}, $\Psi$ has a conformal 1--1 and onto continuous extension to $\overline{\mathbb{D}}$. We have $\Psi'(z) \neq 0$, for $|z|<1$.
 Given that the constant function $ g(z)\equiv 1\in E^2(\mathbb D)$, then by \cite[corollary 10.1]{Duren:1970} we have that $\sqrt{\Psi'}\in H^2(\mathbb D)\subset H^1(\mathbb D)$. Further, $\sqrt{\Psi'}$ is nonzero almost everywhere on $\partial \mathbb D$ (by \cite[Theorem 10.3]{Duren:1970}). 
Therefore we may write $\Psi'(w)=\sqrt{\Psi'(w)} \sqrt{\Psi'(w)}$, for almost every $|w|=1$.
Now, applying the change of variable formula $\zeta:=\Psi(w)$ to $\tau(n)$, we find
\begin{equation}
\begin{split}
\tau(n)&=\int_{\partial \mathbb{D}} {\sqrt{\Phi'(\Psi(w))} (f \circ \Psi)(w) w^n \Psi'(w) \,\dd w} \\
&=\int_{\partial \mathbb{D}} {\sqrt{\Phi'(\Psi(w)) \Psi'(w)} (f \circ \Psi)(w) [\Phi \circ \Psi]^n (w) \sqrt{\Psi'(w)} \,\dd w} \\
&=\int_{\partial \mathbb{D}} {[(f \cdot \Phi^n) \circ \Psi](w) \sqrt{\Psi'(w)} \,\dd w},
\end{split}
\end{equation}
where $\cdot$ denotes the pointwise product, and we have used that $(\Phi \circ \Psi)'(w)=1$, $\forall w$, because $(\Phi \circ \Psi)(w)=w$.

Following Duren \cite{Duren:1970}, we have
\begin{equation}
f \cdot \Phi^n \in E^2(D) \quad \Leftrightarrow \quad [f \cdot \Phi^n] \circ \Psi \cdot \sqrt{\Psi'} \in H^2(\mathbb{D})=E^2(\mathbb{D}),
\label{3.11}
\end{equation}
in which case $\tau(n)=0$, for $n \in \mathbb{Z}^{+}_0$, by Cauchy's theorem for $E^1(\mathbb{D})$ (and the fact that $E^2(\mathbb{D}) \subset E^1(\mathbb{D})$). Therefore, the question whether $\tau(n)=0$ is reduced to whether the following statement is true:
\begin{equation}
f \in E^2(\mathbb{D}) \quad \text{and} \quad \Phi: D \rightarrow \mathbb{D} ~ \text{conformal} \qquad \Rightarrow \qquad  f \cdot \Phi^n \in E^2(D). \label{statement}
\end{equation}

\noindent
To this end, we claim that $\Phi$ extends to $\Phi \in C(\overline{D}) \cap \vartheta(D)$. To see this, by Caratheodory's theorem \cite{Duren:1970}, $\Psi$ extends to a homeomorphism $\Psi:\overline{\mathbb{D}}\mapsto \overline{D}$ (more precisely, $\Psi \in \vartheta(\mathbb{D}) \cap C(\overline{\mathbb{D}})$, $\Psi: \partial \mathbb{D} \mapsto \partial D$ homeomorphism, $\Psi |_{\mathbb{D}}=\psi$). Hence $\Psi^{-1}:\overline{D} \mapsto \overline{\mathbb{D}}$ is a homeomorphism and we may take $\Phi \coloneqq \Psi^{-1}$. The claim is proved.

Now it is well known that (with $\vartheta(D)$ denoting the functions analytic in $D$),
\begin{equation}
E^2(D) \cdot \{ \vartheta(D) \cap C(\overline{D}) \} \subset E^2(D),
\end{equation}
and since $\Phi \in \vartheta(D) \cap C(\overline{D})$, the above implies that the same is true for $\Phi^{n}$, $\forall n \in \mathbb{Z}^{+}_{0}$. It follows that $f \cdot \Phi^n \in E^2(D) \cdot \{ \vartheta(D) \cap C(\overline{D})\} \subset E^2(D)$. Therefore we have shown that the statement \eqref{statement} holds. In order to prove \eqref{GR}, one substitutes the definition of $\rho(k)$ given in \eqref{Transform pair domain2} into the formula for $f(z)$ given in \eqref{Transform pair domain1}. Then \eqref{GR} follows by switching the order of integration and applying Cauchy theorem. See Appendix \ref{Appendix B}. The lemma is proved.
\end{proof}

\section{A transform pair for punctured domains \label{sectiontransformpunctured}}

As a first step to derive a transform pair for a punctured domain $D^{*}$, we consider the conformal mapping $\Phi:D^{*} \mapsto \mathbb{D}^{*}$ which maps the punctured domain $D^{*}$ to the punctured disc $\mathbb{D}^{*}\coloneqq \mathbb{D} \setminus \{0\}$, and we write $\eta=\Phi(\zeta)$ and $w=\Phi(z)$, for $\zeta \in \partial D^{*}$, $z \in D^{*}$.

The function $f$ is in $E^2(D^{*})$ and we assume, without loss of generality, that it has a first-order pole at $\Phi^{-1}(0)$. For the definition of the Hardy space for $D^{*}$ see \cite{GallagherGuptaLanzaniVivas:2021}. From \eqref{E:Szego-punctured-disc}, the Szeg\"o kernel for $\mathbb{D}^{*}$ with a first-order pole at $0$ is given by
\begin{equation}
S_{\mathbb{D^*}}(\eta, w) = \frac{1}{2\pi}\, \frac{1}{\eta \overline{w} (1-\eta\overline{w})}, \quad |\eta|=1, \quad 0<|w|<1.
\end{equation}
Therefore, we can write
\begin{equation}
\overline{S_{\mathbb{D^*}}(\eta, w)} = \frac{1}{2\pi}\, \frac{1}{\overline{\eta} w (1-\overline{\eta}w)} = \frac{1}{2\pi}\, \frac{\eta^2}{w (\eta-w)}.
\end{equation}
Since for $D_2=\mathbb{D}$ we have $T_{D_2}(\Phi(\zeta))=\ci \Phi(\zeta)$, expression \eqref{E:Szego-final} can be written as
\begin{equation}\label{Szego dsigma punctured domain}
\overline{S_{D^{*}}(\zeta,z)}\,\dd \sigma(\zeta) =\sqrt{\Phi'(z)}~\frac{1}{2\pi}\, \frac{\Phi(\zeta)^2}{\Phi(z)(\Phi(\zeta)-\Phi(z))} {\sqrt{\Phi'(\zeta)}} (-\ci \overline{\Phi(\zeta)}) \dd \zeta.
\end{equation}
Substitution of \eqref{Szego dsigma punctured domain} into \eqref{E:SF} gives
\begin{equation}\label{f punctured disc 1}
\begin{split}
f(z)&= \int_{\partial D^{*}} {f(\zeta) \overline{S_{D^{*}}(\zeta,z)} \,\dd \sigma(\zeta)} \\
&=\frac{1}{2\pi \ci} \frac{\sqrt{\Phi'(z)}}{\Phi(z)} \int_{\partial D^{*}} {f(\zeta) \left[ \frac{\Phi(\zeta) \sqrt{\Phi'(\zeta)}}{\Phi(\zeta)-\Phi(z)} \right] \,\dd \zeta}, \qquad \text{for } z \in D^{*}.
\end{split}
\end{equation}
Combining \eqref{f punctured disc 1} with \eqref{E:spectral-general}, we obtain a {\bf transform pair for the punctured domain $D^{*}$ with a first-order pole at $\Phi^{-1}(0)$}, namely:
\begin{equation}\label{Transform pair punctured domain1}
f(z)=\displaystyle \frac{1}{2\pi \ci} \frac{\sqrt{\Phi'(z)}}{\Phi(z)} \left[ \int_{L_{1}}{\frac{\rho(k)}{1-\ee^{2\pi\ci k}}\Phi(z)^{k} \,\dd k}+\int_{L_{2}} {\rho(k) \Phi(z)^{k} \dd k} +\int_{L_{3}}{\frac{\rho(k) \ee^{2\pi\ci k}}{1-\ee^{2\pi\ci k}}\Phi(z)^{k} \,\dd k}\right], \quad z \in D^{*}, 
\end{equation}
where
\begin{equation}\label{Transform pair punctured domain2}
\rho(k)\coloneqq \displaystyle \int_{\partial D^{*}} {\frac{\sqrt{\Phi'(\zeta)}}{\Phi(\zeta)^{k}}f(\zeta) \,\dd \zeta}, \quad k \in \mathbb{C}.
\end{equation}

\noindent
One can show that the global relation is again
\begin{equation}\label{GR punctured}
\rho(k)=0, \quad k \in \mathbb{Z}^{-}.
\end{equation}

\noindent
[Note that \eqref{Transform pair punctured domain1}--\eqref{Transform pair punctured domain2} are different from \eqref{Transform pair domain1}--\eqref{Transform pair domain2} by an extra factor $1/\Phi(z)$ in the representation of $f(z)$ in \eqref{Transform pair punctured domain1}, and an extra factor $1/\Phi(\zeta)$ in the integrand of $\rho(k)$ in \eqref{Transform pair domain2}.]

\section{A boundary value problem on an elliptical domain \label{sectionapplication1}}
In this section, we present a generalization of the mixed boundary value problem considered by Shepherd \cite{Shepherd1937} on the disc to a problem posed on an elliptical domain. The mixed boundary value problem on an elliptical domain was also analysed in our previous study \cite{HLLL:2024} and here it is used for verification of the transform method developed in \S\,\ref{sectiontransform}.

Consider an elliptical domain $D$, namely
\begin{equation}\label{ellipse eqn}
D \coloneqq \left\{z \in \mathbb C, x=\text{Re}[z], y=\text{Im}[z] \left|  \frac{x^2}{a^2}+\frac{y^2}{b^2}<1 \right.\right\},
\end{equation}
whose boundary is
\begin{equation}
\partial D=\left\{z \in \mathbb C, x=\text{Re}[z], y=\text{Im}[z] \left| \frac{x^2}{a^2}+\frac{y^2}{b^2}=1 \right. \right\}.
\end{equation}
We consider the following incomplete Dirichlet boundary value problem for analytic functions
\begin{equation}\label{E:BVP-ell-good}
\begin{cases}
\overline{\partial}f (z) = 0,~~~~~~~~~~~~~~~~~~~~ z \in D, \\
\text{Re} f(\zeta)= \text{Re}\hspace{1mm} \overline{\zeta}^m, \qquad \zeta \in C_1=\{ \zeta\in \partial D \text{ and } \Real \zeta>0 \}, \\
\text{Im} f(\zeta)=\text{Im} \hspace{1mm} \overline{\zeta}^m , \qquad \zeta \in C_2=\{ \zeta\in \partial D \text{ and }  \Real \zeta<0 \},
\end{cases}
\end{equation}
for a given $m\in\mathbb N$. Recall that
\begin{equation}
\overline{\partial}=\frac{1}{2} \left(\frac{\partial }{{\partial x}} +\ci \frac{\partial }{{\partial y}} \right).
\end{equation}

We parametrize the ellipse using polar coordinates $(l,\theta)$, with $0\leq \theta\leq 2\pi$ and
\begin{equation} 
l(\theta)=\frac{ab}{\sqrt{(b\cos\theta)^2+(a\sin\theta)^2}},
\end{equation}
and write       
\begin{equation}\label{parametrization_polar}
\zeta(\theta) \coloneqq l(\theta)\ee^{\ci \theta} \in \partial D.
\end{equation}

\noindent
The boundary value problem \eqref{E:BVP-ell-good} can be written in terms of polar coordinates as
\begin{equation}\label{E:BVP-ell-good_polar}
\begin{cases}
\overline{\partial}f (z) = 0,~~~~~~~~~~~~~~~~~~~~~~~~~~~~~~~~z \in D, \\
\Real f= [l(\theta)]^m \cos m\theta,~~~ \qquad \theta \in (-\pi/2, \pi/2) \quad (\text{or } C_1), \\
\Imag f=- [l(\theta)]^m \sin m\theta, \qquad \theta\in (\pi/2, 3\pi/2) \quad (\text{or } C_2).
\end{cases}
\end{equation}

On $C_1$, we write   
\begin{equation}\label{E:apriori-1}
f(\theta,l(\theta))=[l(\theta)]^m \cos m\theta+\ci\Big(a_0+\sum_{n\geq 1}\Big[a_n\ee^{2\ci n\theta}+\overline{a_n} \ee^{-2 \ci n\theta}\Big]\Big),  
\end{equation}
where the coefficients $a_0 \in \mathbb{R}$ and $\{a_n \in \mathbb C | n=1,2,\dots \}$ are to be determined.

On $C_2$, we write
\begin{equation}\label{E:apriori-2}
f(\theta,l(\theta))= \Big( b_0+\sum_{n\geq 1}\Big[b_n\ee^{2\ci n\theta}+\overline{b_n} \ee^{-2\ci n\theta}\Big] \Big) -\ci [l(\theta)]^m \sin m\theta,
\end{equation}
where the coefficients $b_0 \in \mathbb{R}$ and $\{b_n \in \mathbb C | n=1,2,\dots \}$ are to be determined.

On substitution of \eqref{E:apriori-1} and \eqref{E:apriori-2} into the global relation \eqref{GR}, we obtain a linear system for the unknown coefficients given by
\begin{equation}\label{ellipse:system}
\begin{split}
&a_0\A(0,k)+\sum_{n\geq 1}\Big(a_n\A(2n,k)+\overline{a_n} \A(-2n,k)\Big) \\
&~~~~~~+b_0\B(0,k)+\sum_{n\geq 1}\Big(b_n \B(2n,k)+\overline{b_n} \B(-2n,k)\Big)=r(k), \quad k \in \mathbb{Z}^{-}, \\
&a_0\overline{\A(0,k)}+\sum_{n\geq 1}\Big(\overline{a_n} \overline{\A(2n,k)}+ a_n\overline{\A(-2n,k)}\Big) \\
&~~~~~~+b_0\overline{\B(0,k)} +\sum_{n\geq 1}\Big(\overline{b_n} \overline{\B(2n,k)}+ b_n\overline{\B(-2n,k)}\Big)=\overline{r(k)}, \quad k \in \mathbb{Z}^{-},
\end{split}
\end{equation}
where
\begin{equation}\label{ABellipse}
\A(n,k)= \ci \int_{-\pi/2}^{\pi/2} \frac{\sqrt{\Phi'(\zeta(\theta))}}{\Phi(\zeta(\theta))^{k+1}} \ee^{\ci n \theta} \zeta'(\theta) \,\dd\theta, \qquad \B(n,k)= \int_{\pi/2}^{3\pi/2} \frac{\sqrt{\Phi'(\zeta(\theta))}}{\Phi(\zeta(\theta))^{k+1}} \ee^{\ci n \theta} \zeta'(\theta) \,\dd\theta,
\end{equation}
and the function $r(k)$ is defined by
\begin{equation}\label{rellipse}
r(k)=-\int_{-\pi/2}^{\pi/2} [l(\theta)]^m \cos (m\theta) \frac{\sqrt{\Phi'(\zeta(\theta))}}{\Phi(\zeta(\theta))^{k+1}} \zeta'(\theta)\,\dd\theta
+ \ci \int_{\pi/2}^{3\pi/2} [l(\theta)]^m \sin (m\theta) \frac{\sqrt{\Phi'(\zeta(\theta))}}{\Phi(\zeta(\theta))^{k+1}} \zeta'(\theta)\,\dd\theta.
\end{equation}
The conformal mapping $\Phi$ is defined in Appendix \ref{Appendix A}; it maps the interior of the ellipse $D$ to the unit disc $\mathbb{D}$.

\begin{figure*}[t!]
\hspace*{-1.2cm}
    \centering
    \begin{subfigure}[t]{0.5\textwidth}
        \centering
        \[\includegraphics[scale=0.17] {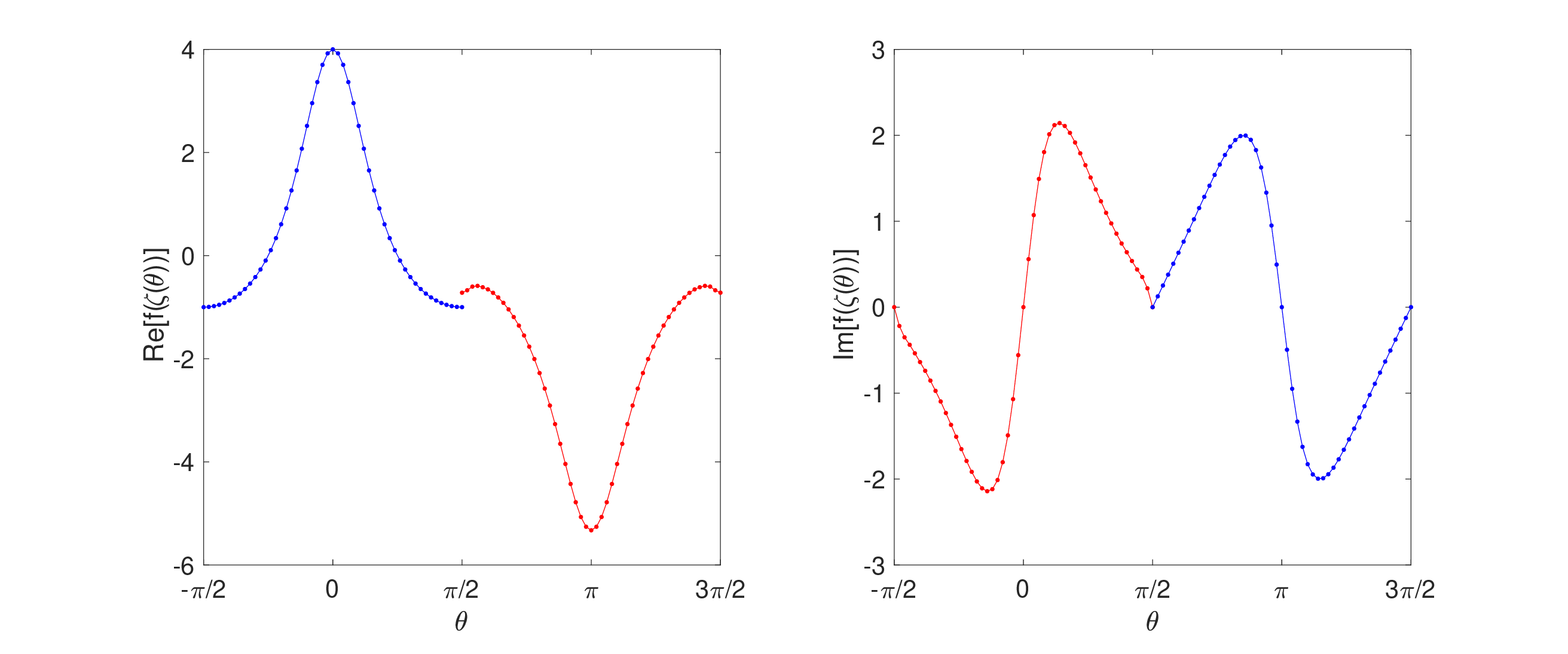} \]
\caption{Parameters: $a=2$, $b=1$.}
\label{Fig:3a}
    \end{subfigure}%
~~   
    \begin{subfigure}[t]{0.5\textwidth}
        \centering
  \[\includegraphics[scale=0.17] {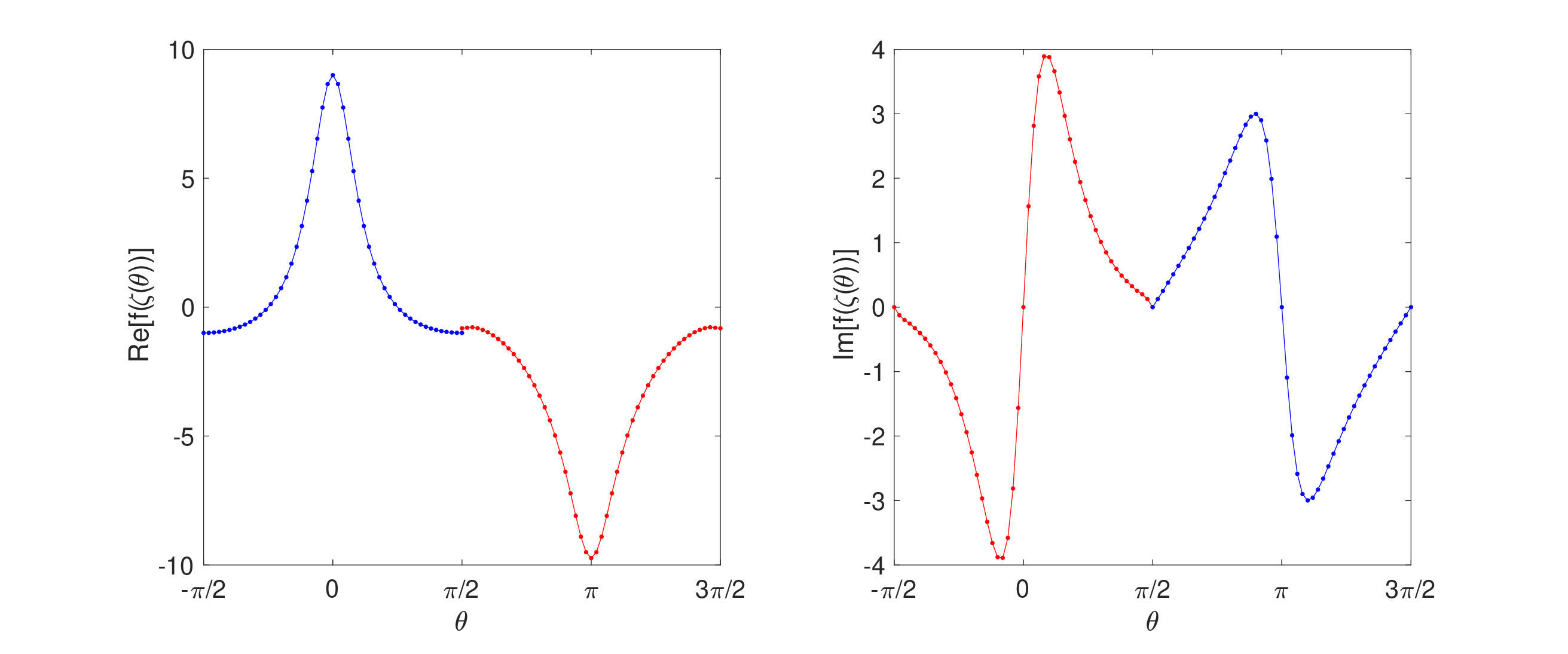} \]
  \caption{Parameters: $a=3$, $b=1$.}
  \label{Fig:3b}
    \end{subfigure}
    \caption{Real and imaginary parts of $f(\zeta(\theta))$, $ \theta \in [-\pi/2, 3\pi/2]$ along the boundary of the ellipse, for $m=2$: numerical solutions \cite{HLLL:2024} (solid lines) and solutions computed via the transform method presented here (dots).}
    \label{fig: ellipse examples}
\end{figure*}

To proceed, the sums \eqref{E:apriori-1} and \eqref{E:apriori-2} are truncated to include only terms up to $n=N$ and we formulate a linear system for the $4N+2$ unknown coefficients. The linear system comprises the conditions \eqref{ellipse:system} evaluated at points $k=-N_k,\dots,-2,-1$, for $N_k>2N+1$, forming an overdetermined linear system of $2N_k$ equations. The unknown coefficients are found by a least-squares algorithm. We found that the coefficients $\{a_n, b_n | n=0,\dots,N \}$ decay quickly and, therefore, we choose the truncation parameter to be $N=16$. Once the coefficients $\{a_n, b_n | n=0,\dots,N \}$ are found, the spectral function $\rho(k)$ and $f(z)$ can be computed via the transform pair \eqref{Transform pair domain1}--\eqref{Transform pair domain2}. Our results here were verified against the numerical solutions that were obtained by means of the transform-based method devised in our previous study \cite{HLLL:2024}.

Figure \ref{fig: ellipse examples} shows the real and imaginary parts of $f(\zeta(\theta))$, $ \theta \in [-\pi/2, 3\pi/2]$ along the boundary of the ellipse, for $m=2$ and different parameter choices for $a$ and $b$, computed via the numerical scheme of \cite{HLLL:2024} and the solutions obtained via the transform method presented here. The numerical solutions via the approaches agree up to 8 decimal digits for similar parameter choices. We observe that the two solutions depicted in Fig.~\ref{fig: ellipse examples} display discontinuities in the real part of $f$ for $\theta=\pi/2$ and $\theta=3\pi/2$ (also corresponding to $\theta=-\pi/2$), which are the values where the boundary conditions change type. Such discontinuities are not surprising given that the solution is in $E^1(D)$.

\section{Application in fluid dynamics: A point vortex inside a concave quadrilateral \label{sectionapplication2}}

In this section, we present an application of the new transform pairs to a problem in fluid dynamics, in particular within the framework of a two-dimensional, inviscid, incompressible, and irrotational (except for point vortices) steady flow. We consider a point vortex (a vortex with infinite vorticity concentrated at a point) in the interior of a concave quadrilateral and the aim is to find the resulting fluid flow satisfying the impermeability boundary condition.

\subsection{Governing equations}

On a given domain $D$, we consider a two-dimensional inviscid, incompressible (the velocity field ${\bf u}=(u,v)$ satisfies the continuity condition $\nabla \cdot {\bf u} = 0$), and irrotational (the vorticity $\omega = \nabla \times {\bf u}$ is zero) steady flow. Since the flow is irrotational, there exists a velocity potential $\phi(x,y)$ such that:
\begin{equation}
u=\frac{\partial \phi}{\partial x}, \qquad v=\frac{\partial \phi}{\partial y}.
\end{equation}
The incompressibility condition then implies that $\phi$ satisfies Laplace's equation:
\begin{equation}
\nabla^2 \phi = 0.
\end{equation}
A stream function $\psi(x,y)$ can also be introduced such that
\begin{equation}
u=\frac{\partial \psi}{\partial y}, \qquad v=-\frac{\partial \psi}{\partial x}.
\end{equation}
The stream function $\psi$ also satisfies Laplace's equation:
\begin{equation}
\nabla^2 \psi = 0.
\end{equation}

\noindent
It is natural to introduce a complex potential function $h(z)$, where
\begin{equation}\label{potential h}
h(z)=\phi(x,y) + \ci \psi(x,y), \quad \text{with } z=x+\ci y.
\end{equation}
The complex potential function $h(z)$ is analytic in $D$ since $\phi(x,y)$ and $\psi(x,y)$ are harmonic conjugates. Note that although $h(z)$ is analytic in the fluid region, isolated singularities are permitted to model flows of interest.

\begin{figure}
\centering
\includegraphics[width=0.9\textwidth]{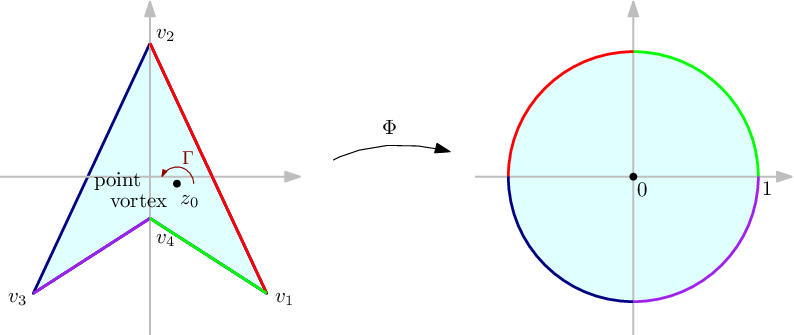}
\caption{Schematic of the configuration: a point vortex at $z_0$ in the interior of the four-sided polygon $D$ (left) and the mapped domain $\mathbb{D}$ (right) via the conformal mapping $\Phi$.}
\label{vortex 4sided}
\end{figure}

\subsection{Problem formulation and solution scheme \label{quadrilateralvortex1st}}

Consider a point vortex with circulation $\Gamma$ at point $z_0$ in the interior of the concave quadrilateral $D$ with vertices $v_j$, $j=1$, $2$, $3$, $4$, as shown in Fig.~\ref{vortex 4sided}. To begin, we introduce the complex potential function $h(z)$ as in \eqref{potential h} and write
\begin{equation}\label{cp}
h(z)=f_s(z)+f(z),
\end{equation}
where
\begin{equation}
f_s(z)=\frac{\Gamma}{2 \pi \ci} \log(z-z_0).
\end{equation}
The function $f(z)$ is analytic in $D$ and will be found using the transform method developed in \S\,\ref{sectiontransform}.

We impose an impermeability condition on the boundary of $D$ which means that the normal velocity component must vanish on each segment of the boundary. This is equivalent to requiring that the stream function $\psi$ be constant along the boundary; without loss of generality we assume that this constant is zero. The impermeability condition can be expressed in terms of the complex potential as
\begin{equation}\label{BCellipse}
\psi=\text{Im}[h(\zeta)]=0, \quad \zeta \in \partial D.
\end{equation}

Substitution of \eqref{cp} into \eqref{BCellipse} gives
\begin{equation}
\text{Im}[f(\zeta)]=-\text{Im}[f_s(\zeta)].
\end{equation}
We represent $f(\zeta)$ on the four sides $S_j$, $j=1$, $2$, $3$, $4$ of $\partial D$ using Chebyshev expansions:
\begin{equation}\label{function rep1}
f(\zeta_j(s)) = \sum_{n=0}^{\infty} a_{jn} T_n(s) - \ci ~\text{Im}[f_s(\zeta_j(s))], \quad \text{for } s \in [-1,1], \qquad j=1,2,3,4, 
\end{equation}
where $T_{n}(s)=\cos(n \cos^{-1}(s))$ is the $n$th Chebyshev polynomial and $\{a_{jn} \in \mathbb{R} | n=0,1,2,\dots \}$ are to be determined.
The sides $S_j$, $j=1,2,3,4$ are parametrized by
\begin{equation}\label{quadrilateral_sides_param}
\zeta_j(s) = \frac{1-s}{2} v_{j} + \frac{1+s}{2} v_{j+1}, \quad s \in [-1,1],
\end{equation}
respectively, with vertices $v_j$, $j=1$, $2$, $3$, $4$ (and $v_5=v_1$), as shown in Fig.~\ref{vortex 4sided}.

The infinite sums in \eqref{function rep1} are truncated to include terms up to $n=N$. On substitution of \eqref{function rep1} into the global relation \eqref{GR}, we obtain (after some algebra and rearrangement) a linear system for the $4N+4$ unknown coefficients $\{a_{jn} \in \mathbb{R}| n=0,\dots,N\}$, for $j=1$, $2$, $3$, $4$. The interior of $D$ is mapped to the unit disc $\mathbb{D}$ via a Schwarz--Christoffel mapping $\Phi$ computed numerically \cite{SCbook}. The linear system is given by

\begin{equation}
\sum_{j=1}^{4} \sum_{n=0}^{N} a_{jn} P_{j}(n,k)=R(k), \quad \text{for } k \in \mathbb{Z}^{-},
\label{linearsystempointsingularity}
\end{equation}
where
\begin{equation}
P_{j}(n,k)= \int_{-1}^{1} {T_n(s) \frac{\sqrt{\Phi'(\zeta_j(s))}}{\Phi(\zeta_j(s))^{k+1}} \zeta_j'(s) \,\dd s}
\end{equation} 
and
\begin{equation}
R(k) = \ci \sum_{j=1}^{4} \int_{-1}^{1} {\text{Im}[f_s(\zeta_j(s))] \frac{\sqrt{\Phi'(\zeta_j(s))}}{\Phi(\zeta_j(s))^{k+1}} \zeta_{j}'(s) \,\dd s}.
\end{equation}
The linear system \eqref{linearsystempointsingularity} is evaluated at points $k=-N_k,\dots,-2,-1$, for $N_k>4N+4$, which are used to form an overdetermined linear system. The unknown coefficients are found by a least-squares algorithm. We found that the coefficients $\{a_{jn} \in \mathbb{R}| n=0,\dots,N \}$, for $j=1, 2, 3, 4$, decay quickly and, therefore, we choose the truncation parameter to be $N=16$. Once the coefficients are found, the spectral functions and $f(z)$ can be computed via the transform pair \eqref{Transform pair domain1}--\eqref{Transform pair domain2}. Our results were verified against the exact solution for the complex potential function given by
\begin{equation}
f(w)=\frac{\Gamma}{2\pi\ci} \log \left(\frac{w-w_0}{w-1/\overline{w_0}}\right), \quad \text{with }\ w_0=\Phi(z_0),
\end{equation}
where $w=\Phi(z)$. Figure \ref{polygonstreamlines} shows the streamline pattern for a point vortex of strength $\Gamma=1$ at points in the interior of the quadrilateral $D$.

\begin{figure}
\hspace*{-0.8cm}
\centering
\includegraphics[width=1\textwidth]{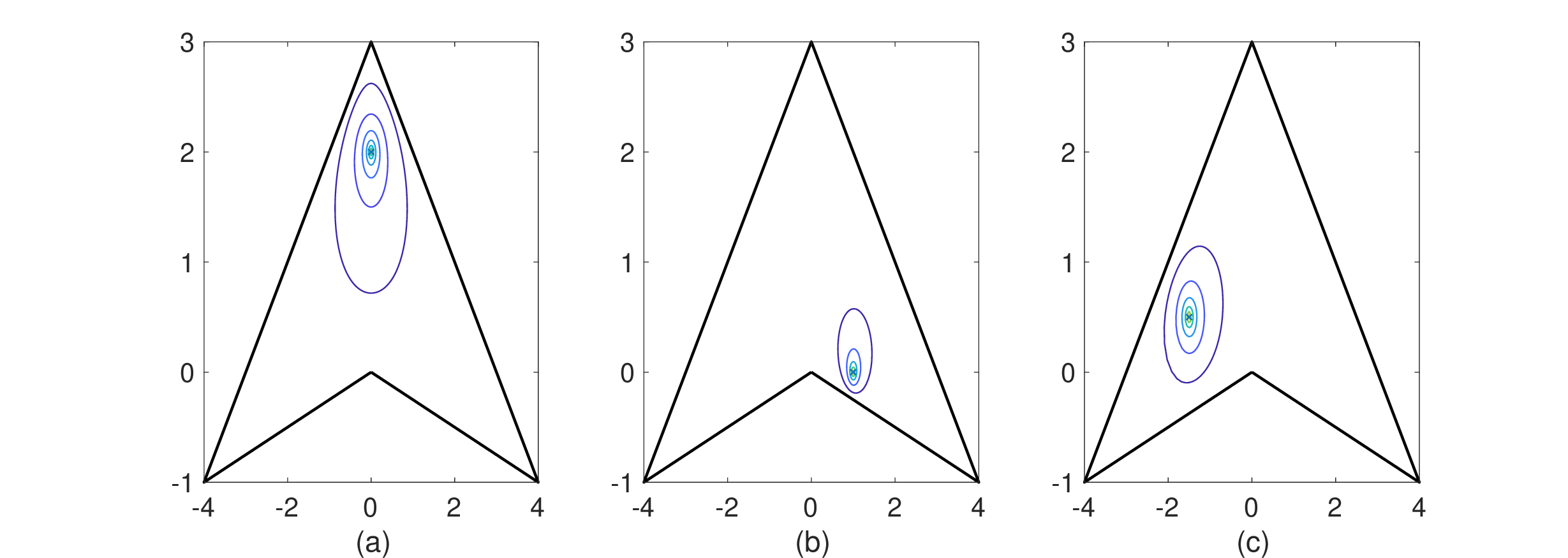}
\caption{Streamlines for a point vortex of strength $\Gamma=1$ at point $z_0$ in the interior of the four-sided polygon $D$ with vertices $\{a-\ci c,\ci b,-a-\ci c,0 \}$, for $(a,b,c)=(4,3,1)$. Results are shown for the point vortex at (a) $z_0=2\ci$, (b) $z_0=1$ and (c) $z_0=-1.5+0.5\ci$.}
\label{polygonstreamlines}
\end{figure}

\begin{figure}
\hspace*{-0.8cm}
\centering
\includegraphics[width=1\textwidth]{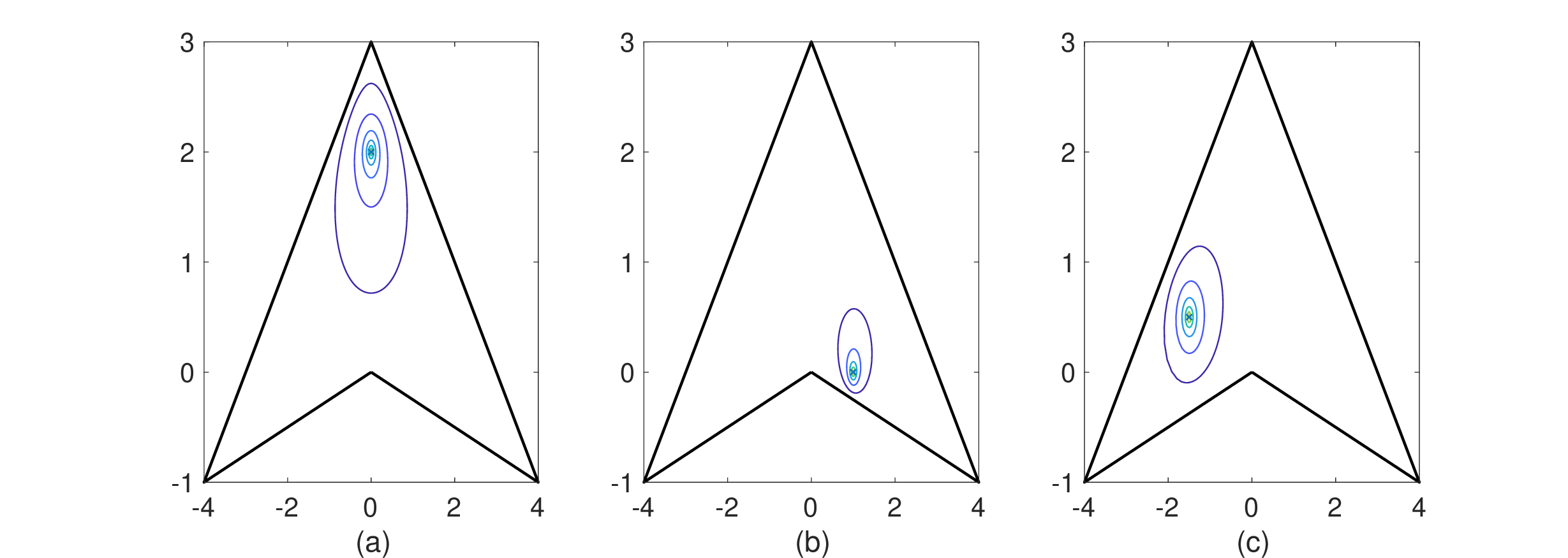}
\caption{Streamlines for a point vortex of strength $\Gamma=1$ at point $z_0$ in the interior of the four-sided polygon $D$ with vertices $\{a-\ci c,\ci b,-a-\ci c,0 \}$, for $(a,b,c)=(4,3,1)$. Results are shown for the point vortex at (a) $z_0=2\ci$, (b) $z_0=1$ and (c) $z_0=-1.5+0.5\ci$.}
\label{polygonstreamlines}
\end{figure}

\subsection{Alternative approach via mapping to the punctured disc \label{quadrilateralvortex2nd}}

Can the transform pair for punctured domains provide a more effective approach to the solution of the point vortex problem on $D$ that was considered in \S\,\ref{quadrilateralvortex1st}? To this end, we reformulate the original problem on $D$ as a problem on the punctured domain $D^*$, and correspondingly use the transform pair presented in \S\,\ref{sectiontransformpunctured}. Specifically, we consider the punctured concave quadrilateral $D^*$ and its corresponding Szeg\H o kernel, and we use a conformal mapping $\Phi$ from $D^*$ to the punctured disc $\mathbb{D}^{*}$. The puncture is located at point $z_0=\Phi^{-1}(0)$ and assume we have a point vortex with circulation $\Gamma$ at point $z_0$.

We introduce the function $w(z)$:
\begin{equation}\label{alternative w(z) pv}
w(z)=w_s(z) + \frac{\Gamma}{2\pi \ci} f(z),
\end{equation}
where
\begin{equation}
w_s(z) = \frac{\Gamma}{2 \pi \ci} \frac{1}{z-z_0},
\end{equation}
where $f(z)$ is analytic in $D^*$ (as opposed to the function $f(z)$ analytic in $D$ that was considered in \S\,\ref{quadrilateralvortex1st}), and will be found using the transform method presented in \S\,\ref{sectiontransformpunctured}. The function $w_s(z)$ is analytic in $D^{*}$. Note that the function $w(z)$ represents the complex velocity $w=u-\ci v$.

The impermeability boundary condition on the quadrilateral can be written as
\begin{equation}
\text{Im}\left[{\bf u} \cdot {\bf n} \right]=0,
\end{equation}
where ${\bf n}=(n_x,n_y)$ is the unit normal to the boundary. This can be expressed as
\begin{equation}\label{bc alternative pv}
\text{Im}\left[w(\zeta) (-\ci) \zeta'(s) \right]=0.
\end{equation}
Substitution of \eqref{alternative w(z) pv} into \eqref{bc alternative pv} gives
\begin{equation}
\text{Im}\left[\frac{\Gamma}{2 \pi} f(\zeta) \zeta'(s) \right]= \text{Im} [w_s(\zeta) (-\ci) \zeta'(s) ].
\end{equation}

Proceeding in a fashion analogous to \eqref{E:apriori-1}, \eqref{E:apriori-2} and \eqref{function rep1}, we write
\begin{equation}\label{function rep1 alternative}
\frac{\Gamma}{2\pi} f(\zeta_j(s)) \zeta_j'(s) = \sum_{n=0}^{\infty} b_{jn} T_n(s) + \ci ~\text{Im}\left[w_s(\zeta_j(s)) (-\ci) \zeta_j'(s) \right], \quad s \in [-1,1],
\end{equation}
on the four sides $S_j$, $j=1$, $2$, $3$, $4$ of $\partial D^{*}$, where the coefficients $\{b_{jn} \in \mathbb{R} | n=0,1,2,\dots \}$ are to be determined.

On substitution of \eqref{function rep1 alternative} into the global relation \eqref{GR punctured}, we obtain (after some algebra and rearrangement) a linear system for the $4N+4$ unknown coefficients $\{b_{jn} \in \mathbb{R} | n=0,\dots,N\}$, for $j=1$, $2$, $3$, $4$. The infinite sums are truncated to include terms up to $n=N$. The linear system is given by

\begin{equation}
\sum_{j=1}^{4} \sum_{n=0}^{N} b_{jn} U_j(n,k) =V(k), \quad \text{for } k \in \mathbb{Z}^{-},
\label{linearsystempointsingularity alternative}
\end{equation}
where
\begin{equation}
U_j(n,k)=\int_{-1}^{1} {T_n(s) \frac{\sqrt{\Phi'(\zeta_j(s))}}{\Phi(\zeta_j(s))^{k}} \,\dd s}
\end{equation}
and
\begin{equation}
V(k) = - \ci \sum_{j=1}^{4} \int_{-1}^{1} {\text{Im}[w_s(\zeta_j(s)) (-\ci) \zeta'_j(s)] \frac{\sqrt{\Phi'(\zeta_j(s))}}{\Phi(\zeta_j(s))^{k}} \,\dd s}.
\end{equation}
The linear system \eqref{linearsystempointsingularity alternative}  is evaluated at points $k=-N_k,\dots,-2,-1$, for $N_k>4N+4$, which are used to form an overdetermined linear system. The unknown coefficients are found by a least-squares algorithm. We found that the coefficients $\{b_{jn} \in \mathbb{R}| n=0,\dots,N \}$, for $j=1$, $2$, $3$, $4$, decay quickly and, therefore, we choose the truncation parameter to be $N=16$. Once the coefficients are found, the spectral functions and $f(z)$ can be computed via the transform pair \eqref{Transform pair punctured domain1}--\eqref{Transform pair punctured domain2}. This alternative approach does not appear to present any advantages nor disadvantages in comparison with the approach described in \S\,\ref{quadrilateralvortex1st}, namely, the numerical computations agree up to 16 decimal digits for similar parameter choices.

\section{Summary and discussion \label{sectiondiscussion}}

In this study, we have expanded the scope of the UTM to address boundary value problems in Lipschitz and in particular non-convex planar domains, including shapes beyond circular. Building on the work by Crowdy \cite{Crowdy:2015} for circular domains, we developed a novel transform method which is also applicable to punctured domains. We analysed a mixed boundary value problem in an elliptical domain with different boundary conditions imposed on different arcs of the boundary, as well as a fluid dynamics problem inside a non-convex punctured quadrilateral; for the latter problem we have in fact pursued two different approaches in \S\,\ref{quadrilateralvortex1st}--\ref{quadrilateralvortex2nd} and neither appears to be more or less advantageous than the other.

Our work underscores the importance of the Szeg\H o kernel and its connection to the Cauchy kernel in extending the UTM beyond polygonal or circular domains. The Cauchy kernel lacks favourable transformation properties under conformal maps, while the Szeg\H o kernel has a good transformation law. The applicability of this new transform pair relies on the numerical computations afforded by the conformal mapping $\Phi$ from the domain of interest to the unit disc. This provides a framework that broadens the applicability of the UTM to any simply-connected Lipschitz domain, including non-convex domains, with a conformal map to the unit disc, as well as punctured Lipschitz domains.

\begin{figure}
\centering
\includegraphics[width=0.9\textwidth]{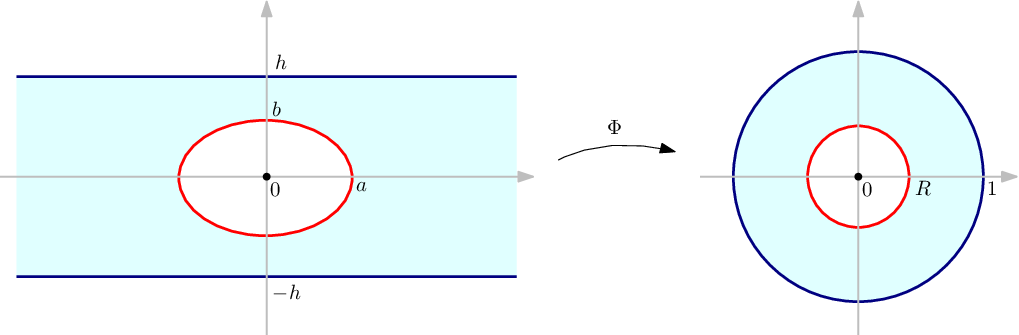}
\caption{Domain $D$ is the region exterior to the ellipse \eqref{ellipse eqn} and interior to the channel $-\infty<x<\infty$, $-h<y<h$, with $b<h$ (left). The mapped domain $\mathbb{A}$ is the annulus with $0<R<1$ (right).}
\label{ellipse in channel}
\end{figure}

The transform method developed in this paper serves as a stepping stone for further research in solving boundary value problems in more complex geometries, including {\it exterior/unbounded domains} and {\it multiply-connected domains}. The following questions naturally arise:
\begin{itemize}
\item {\it Exterior/unbounded domains:} Given a bounded domain $D$, we denote its complement by $D_E=\mathbb{C} \setminus \overline{D}$. Note that $\partial D_{E}= \partial D$ and $\dd\sigma_{D_E}(z)=-\dd \sigma_{D}(z)$. Is it possible to obtain a Szeg\H o projection/kernel for $D_{E}$ via a conformal map $\Phi: \mathbb{C}\setminus \overline{D} \mapsto \mathbb{C} \setminus \overline{\mathbb{D}}$? The Szeg\H o formula
\begin{equation}
f(z)=\int_{\partial D_{E}} {f(\zeta) \overline{S_{D_E}(\zeta,z)} \,\dd \sigma_{D_E}(\zeta)}, \qquad z \in D_{E}=\mathbb{C}\setminus \overline{D}
\end{equation}
is valid for any $f \in E^2(D_E)$.
\item \emph{Multiply-connected domains:} Consider a multiply-connected domain $D$. Without loss of generality, assume that $D$ is a doubly-connected domain (for example, the domain $D$ could be the region exterior to the ellipse \eqref{ellipse eqn} and interior to the channel $-\infty<x<\infty$, $-h<y<h$, with $b<h$); a schematic is shown in Fig.~\ref{ellipse in channel}. Assume we can determine numerically a conformal mapping $\Phi: D \mapsto \mathbb{A}$, where $\mathbb{A}$ is the annulus with $0<R<1$. The Szeg\H o kernel for the annulus is given by \cite{TegtmeyerThomas:1999}:
\begin{equation}
S_{a}(z)=\frac{1}{2\pi} \sum_{n=-\infty}^{\infty} \frac{(z\overline{a})^n}{1+R^{2n+1}}.
\end{equation}
Is it possible to construct a transform pair for the doubly-connected domain $D$? Is it possible to extend this approach to domains with higher-connectivity?
\end{itemize}

Finally, we note that this new approach can be also used to analyse boundary value problems for the biharmonic equation (which will involve solving for two analytic functions; Langlois \cite{Langlois1964}, Luca \& Crowdy \cite{LucaCrowdy:2018}, as well as the complex Helmholtz equation discussed by Hauge \& Crowdy \cite{HaugeCrowdy:2021} and Hulse {\it et al.}~\cite{HLLL:2024}.

\enlargethispage{20pt}

\section*{Acknowledgments }The authors would like to thank the Isaac Newton Institute for Mathematical Sciences, Cambridge, for support and hospitality during the programme {\it Complex analysis: techniques, applications and computations} where part of the work on this paper was initiated, as well as the workshop {\it Complex analysis: techniques, applications and computations - perspectives in 2023}.

\section*{Funding} This work was partly supported by EPSRC grant no EP/R014604/1. L. Lanzani is a member of the INdAM group GNAMPA.

\appendix

\section{Conformal mapping for the ellipse \label{Appendix A}}

Consider the ellipse
\begin{equation}
\frac{x^2}{a^2}+\frac{y^2}{b^2}=1, \label{ellipse}
\end{equation}
with $a>b$ and foci at distance $c=\sqrt{a^2-b^2}>0$.

{\bf Mapping the interior of an ellipse to the disc:} The conformal mapping $\Phi$ from the ellipse \eqref{ellipse} in the $z$-plane to the unit disc $\mathbb{D}$ is given by (Schwarz \cite{Schwarz:1869}):
\begin{equation}\label{ellipse interior map 1}
\Phi(z) = \sqrt{k} ~\text{sn} \left(\frac{2K}{\pi} \sin^{-1}(z),m \right),
\end{equation}
where $\text{sn}$ is the Jacobi sn function and the definition of parameters $k,K$ and $m$ is given in \cite{Schwarz:1869,Szego:1950}.

\section{Switching the order of integration in the transform pairs \label{Appendix B}}

Let $\dd |\zeta|$ and $\dd |k|$ denote the arc-length measures with respect to $\zeta$ and $k$. For the new transform pair, the order of integration may be switched if
\begin{equation}
    \displaystyle\int_{L_{1}}{\displaystyle \int_{\partial D} {\Big|\frac{\sqrt{\Phi'(\zeta)}}{\Phi(\zeta)^{k+1}}}\frac{f(\zeta)}{1-\ee^{2\pi\ci k}}\Phi(z)^{k} \Big|\,\dd |k| \dd |\zeta|}<\infty
\end{equation}
with analogous inequalities also for $L_2$ and $L_3$. Further, we only need to check that the modulus of the integrand decays sufficiently as $|k|\rightarrow \infty$. Hence it suffices to check the integrals above are finite with $L_j$ replaced with $L^+_j:=L_j\cap \{ |k|>r\}$ for $1\leq j\leq 3$.
 Let $\text{Arg}(\zeta)$ denote the principle argument of $\zeta$. The complex exponential is multi-valued, in particular when $|w|=1$, we have that
 \begin{equation}
 w^{2k+1}=\ee^{2 \ci k(\text{Arg}(w)+2\pi m)}, \quad m \in \mathbb Z.
 \end{equation}
 Choose $m=0$ and let $k\in L_1^+$. Recall that $L_1^+=\{ -\ci |k|:|k|\geq r\}$ (see Fig.~\ref{fundamental contour}).  Denote $\theta:=\text{Arg}(w)$. Then 
\begin{equation}
\Big|\frac{1}{w^{2(k+1)}}\Big|=\Big|\frac{1}{\ee^{2\ci k\text{Arg}(w)}}\Big|=\ee^{-2|k|\text{Arg}(w)}=\ee^{-2\theta|k|}.
\end{equation}
Applying the change of variables $\zeta=\Psi(w)$, we find
\begin{equation}
\rho(k)=\int_{\partial D}\frac{\sqrt{\Phi'(\zeta)}\,f(\zeta)}{\Phi(\zeta)^{k+1}}\,\dd \zeta=\int_{\partial \mathbb D} \frac{f(\Psi(w))}{w^{k+1}}\sqrt{\Psi'(w)}\,\dd w.
\end{equation}
Let
\begin{align}
    a\leq \text{Arg}(w)=\theta<b \text{ where } b-a=2\pi.
\end{align}
Holder's inequality yields
\begin{equation}
\begin{split}
\int_{\partial \mathbb D} \Big|\frac{f(\Psi(w))}{w^{k+1}}\sqrt{\Psi'(w)}\Big|\,\dd| w| &\leq\sqrt{ \int_{\partial \mathbb D} \Big|\frac{1}{w^{2(k+1)}}\Big|\dd|w|}\sqrt{\int_{\partial \mathbb D} |f^2(\Psi(w))\Psi'(w)\,|\dd|w|} \\
&\leq M\sqrt{\int_{\partial \mathbb D} \Big|\frac{1}{w^{2(k+1)}}\,\Big|\dd|w|}
 =M\sqrt{\Big[\frac{1}{2|k|} \ee^{-2\theta |k|}\Big]_{\theta=a}^{\theta=b}}.
\end{split}
\end{equation}
We note that $H^2(\mathbb D)\subset H^1(\mathbb D)$ so that $f^2(\Psi(w))\Psi'(w)\in E^1(\mathbb D)=H^1(\mathbb D)$ by \cite[corollary 10.1]{Duren:1970}.
If we let $a=0$ and $b=2\pi$, then we have
\begin{equation}
\oint_{\partial \mathbb D} \Big|\frac{f(\Psi(w))}{w^{k+1}}\sqrt{\Psi'(w)}\Big|\dd| w| \leq M\sqrt{\Big[\frac{1}{2|k|} \ee^{-2\theta |k|}\Big]_{\theta=-\pi}^{\theta=\pi}}= \frac{M}{\sqrt{2|k|}}\frac{\sqrt{|1-\ee^{4\pi |k| }|}}{\ee^{\pi |k| }}.
\end{equation}   
Combining the observations above gives
\begin{equation}\label{Loneplus}
\begin{split}
\int_{L_1^+}\int_{\partial D}  \Big| \frac{\Phi(z)^k}{1-\ee^{2 \ci \pi 
k}}\frac{\sqrt{\Phi'(\zeta)}\,f(\zeta)} {\Phi(\zeta)^{k+1}}\Big| \dd |\zeta| \dd |k| &= \int_{L_1^+}\Bigg|  \frac{\Phi(z)^k}{1-\ee^{2 \ci \pi 
k}}\Bigg|\Bigg[\int_{\partial \mathbb D}  \Big|\frac{\sqrt{\Phi'(\zeta)}\,f(\zeta)} {\Phi(\zeta)^{k+1}}\Big| \dd |\zeta| \Bigg]\dd |k| \\
&\leq \int_{r}^{\infty}\Big|{\frac{\Phi(z)^{k}}{1-\ee^{2\pi  |k|}}\Big[\frac{M}{\sqrt{2|k|}}\frac{\sqrt{|1-\ee^{4\pi |k| }|}}{\ee^{\pi |k| }}\Big]\Big|\dd |k|} \\
&=\int_{r}^{\infty}\frac{M}{\sqrt{2|k|}}\Big|\frac{\sqrt{|1-\ee^{4\pi |k| }|}}{1-\ee^{2\pi  |k|}}\Big| \ee^{-|k|(\pi -\text{arg}(\Phi(z))}\dd |k|.
\end{split}
\end{equation}
For the integral at the end of (\ref{Loneplus}) to be finite, we need
\begin{equation}
-|k|(\pi -\text{arg}(\Phi(z)))< 0.
\end{equation}
Given that $-|k|<0$, this holds when $\text{arg}(\Phi(z))<\pi$. Repeating the work above for $L_3^+$ shows that we may switch the order of integration when $\text{arg}(\Phi(z))>-\pi$. The order of integration on $L_2^+$ may be switch regardless of $\text{arg}(\Phi(z))$. Thus the order of integration may be switched if $-\pi <\text{arg}(\Phi(z))<\pi$. Repeating the work above with $a=0$ and $b=2\pi$ shows that the order of integration may be switched if $0<\arg(\Phi(z))<2\pi$.


\vskip2pc

\bibliography{Szegobib12a} 

\begin{thebibliography}{10}

\bibitem{Bell:1992}
S.~R. Bell.
\newblock {\em The Cauchy transform, potential theory, and conformal mapping}.
\newblock Chapman and Hall/CRC, New York, 2nd edition, 2015.

\bibitem{Conway:1990}
J.~B. Conway.
\newblock {\em A course in functional analysis}.
\newblock Springer-Verlag, New York, 1990.

\bibitem{Crowdy:2015}
D.~G. Crowdy.
\newblock {Fourier--Mellin} transforms for circular domains.
\newblock {\em Comput. Methods Funct. Th.}, 15(4):655 -- 687, 2015.

\bibitem{Crowdy:2015b}
D.~G. Crowdy.
\newblock A transform method for {Laplace's} equation in multiply connected circular domains.
\newblock {\em IMA J. Appl. Math.}, 80:1902--1931, 2015.

\bibitem{CrowdyFokas2004}
D.~G. Crowdy and A.~S. Fokas.
\newblock Explicit integral solutions for the plane elastostatic semi-strip.
\newblock {\em Proc. R. Soc. Lond. A}, 460:1285--1310, 2004.

\bibitem{CrowdyLuca:2014}
D.~G. Crowdy and E.~Luca.
\newblock Solving {Wiener--Hopf} problems without kernel factorization.
\newblock {\em Proc. R. Soc A}, 470:20140304, 2014.

\bibitem{DavisFornberg2014}
C.~Davis and B.~Fornberg.
\newblock A spectrally accurate numerical implementation of the fokas transform method for helmholtz-type pdes.
\newblock {\em Complex. Var. Elliptic}, 59:564--577, 2014.

\bibitem{DimakosFokas2015}
M.~Dimakos and A.~S. Fokas.
\newblock The poisson and the biharmonic equations in the interior of a convex polygon.
\newblock {\em Stud. Appl. Math.}, 134:456--498, 2015.

\bibitem{SCbook}
T.~A. Driscoll and L.~N. Trefethen.
\newblock {\em Schwarz--Christoffel Mapping}.
\newblock Cambridge University Press, 2002.

\bibitem{Duren:1970}
P.~L. Duren.
\newblock {\em Theory of {$H^p$} spaces}.
\newblock Academic Press, New York, 1970.

\bibitem{Fokas1997}
A.~S. Fokas.
\newblock A unified transform method for solving linear and certain nonlinear pdes.
\newblock {\em Proc. Roy. Soc. London Ser. A}, 453:1411--1443, 1997.

\bibitem{FokasKapaev2003}
A.~S. Fokas and A.~A. Kapaev.
\newblock On a transform method for the laplace equation in a polygon.
\newblock {\em IMA J. Appl. Math.}, 68:355--408, 2003.

\bibitem{GallagherGuptaLanzaniVivas:2021}
A.-K. Gallagher, P.~Gupta, L.~Lanzani, and L.~Vivas.
\newblock Hardy spaces for a class of singular domains.
\newblock {\em Math. Z.}, 299:2171--2197, 2021.

\bibitem{HaugeCrowdy:2021}
J.~C. Hauge and D.~G. Crowdy.
\newblock A new approach to the complex helmholtz equation with application to diffusion wave fields, impedance spectroscopy and unsteady stokes flow.
\newblock {\em IMA J. Appl. Math.}, 86:1287--1326, 2021.

\bibitem{HLLL:2024}
J.~J. Hulse, L.~Lanzani, S.~G. {Llewellyn Smith}, and E.~Luca.
\newblock A transform pair for bounded convex planar domains.
\newblock {\em IMA J. Appl. Math.}, 89:574--597, 2024.

\bibitem{KerzmanStein:1978}
N.~Kerzman and E.~M. Stein.
\newblock The szegő kernel in terms of cauchy--fantappi\'e kernels.
\newblock {\em Duke Mathematical Journal}, 45(2):197--224, 1978.

\bibitem{Langlois1964}
W.~E. Langlois.
\newblock {\em Slow viscous flows}.
\newblock Macmillan, New York, NY, 1964.

\bibitem{LucaCrowdy:2018}
E.~Luca and D.~G. Crowdy.
\newblock A transform method for the biharmonic equation in multiply connected circular domains.
\newblock {\em IMA J. Appl. Math.}, 83:942--976, 2018.

\bibitem{Pommerenke:1992}
C.~Pommerenke.
\newblock {\em Boundary behaviour of conformal maps}.
\newblock Springer-Verlag, Berlin; New York, 1992.

\bibitem{Schwarz:1869}
H.A. Schwarz.
\newblock Ueber einige {A}bbildungsaufgaben.
\newblock {\em J. Reine Angew. Math.}, 70:105--120, 1869.

\bibitem{Shepherd1937}
W.~M. Shepherd.
\newblock On trigonometric series with mixed conditions.
\newblock {\em Proc. Lond. Math. Soc.}, 2(43):366--375, 1937.

\bibitem{SpenceFokas2010}
E.~Spence and A.~S. Fokas.
\newblock A new transform method i: domain-dependent fundamental solutions and integral representations.
\newblock {\em Proc. R. Soc. A}, 466:2259--2281, 2010.

\bibitem{Szego:1950}
G.~{Szegő}.
\newblock Conformal mapping of the interior of an ellipse onto a circle.
\newblock {\em The American Mathematical Monthly}, 57(7):474--478, 1950.

\bibitem{TegtmeyerThomas:1999}
T.~J. Tegtmeyer and A.~D. Thomas.
\newblock The ahlfors map and szegő kernel for an annulus.
\newblock {\em Rocky Mountain J. Math.}, 29(2):709--723, 1999.

\end{thebibliography}
\bibliographystyle{plain}

\end{document}